\documentclass[11pt,twoside]{article}

\usepackage{latexsym, amssymb, amsmath, theorem}
\numberwithin{equation}{section}

\theoremstyle{plain} \theorembodyfont{\itshape}
\newtheorem{theorem}{Theorem}[section]
\newtheorem{proposition}[theorem]{Proposition}

\theorembodyfont{\rmfamily}
\newtheorem{definition}[theorem]{Definition}

\newtheorem{remark}[theorem]{Remark}

\usepackage{amscd,amsfonts}
\usepackage{pst-all}
\usepackage{amsfonts}
\usepackage{epsf,graphicx}
\usepackage{epic}

\usepackage{amssymb, upref}
\usepackage{epsf,subfigure,  verbatim}

 \def\eps{{\epsilon}}
\def\C{{\mathbb C}}

\def\N{{\mathbb N}}

\def\ds{\displaystyle}

\title{Divergent series: past, present, future \dots
\thanks{This work is supported by NSERC in Canada.}}
\author{Christiane ROUSSEAU\\[5pt]
\textit{D\'epartement de math\'ematiques et de statistique and
CRM}
\\ \textit{Universit\'e de Montr\'eal}
\\ \textit{C.P. 6128, Succursale Centre-ville, Montr\'eal (Qu\'e.), H3C 3J7,
Canada}}

\begin{document}
\maketitle

\begin{abstract} The paper presents some  reflections of
the author on divergent series and their role and place in
mathematics over the centuries. The point of view presented here
is limited to differential equations and dynamical systems.
\end{abstract}

\section{Introduction}

In the past, divergent series have played a central role in
mathematics. Mathematicians  like, for instance, Lacroix, Fourier,
Euler, Laplace, etc. have used them extensively. Nowadays, they
play a marginal role in mathematics and are often not mentioned in
the standard curriculum. A significant number of mathematicians and most
students are not aware that they can be of any use.

\bigskip
\begin{center}
Where was the turn?
\end{center}
\bigskip According to \cite{B}, the turn occured at the time of Cauchy and Abel, when the
need was felt for constructing analysis on absolute rigor.

\bigskip So let us go back and see what Cauchy and Abel
have said of divergent series.
\bigskip

Cauchy (Preface of \lq\lq Analyse math\'ematique", 1821): {\it
\lq\lq J'ai \'et\'e forc\'e d'admettre diverses propositions qui
para\^\i tront peut-\^etre un peu dures. Par exemple qu'une
s\'erie divergente n'a pas de somme\ldots "} (\lq\lq I have been
forced to admit some propositions which will seem, perhaps, hard
to accept. For instance that a divergent series has no sum
\ldots")

\bigskip

Cauchy made one exception: he justified rigorously the use of the
divergent Stirling series to calculate the $\Gamma$-function. We
will explain below the kind of argument he made when looking at
the example of the Euler differential equation.

\bigskip
Abel (Letter to  Holmboe, January 16 1826): {\it \lq\lq Les
s\'eries divergentes sont en g\'en\'eral quelque chose de bien
fatal et c'est une honte qu'on ose y fonder aucune
d\'emonstration. On peut d\'emontrer tout ce qu'on veut en les
employant, et ce sont elles qui ont fait tant de malheurs et qui
ont enfant\'e tant de paradoxes. \ldots Enfin mes yeux se sont
dessill\'es d'une mani\`ere frappante, car \`a l'exception des cas
les plus simples, par exemple les s\'eries g\'eom\'etriques, il ne
se trouve dans les math\'ematiques presque aucune s\'erie infinie
dont la somme soit d\'etermin\'ee d'une mani\`ere rigoureuse,
c'est-\`a-dire que la partie la plus essentielle des
math\'ematiques est sans fondement. Pour la plus grande partie les
r\'esultats sont justes il est vrai, mais c'est l\`a une chose
bien \'etrange. Je m'occupe \`a en chercher la raison, probl\`eme
tr\`es int\'eressant."} (\lq\lq  Divergent series are, in general,
something terrible and it is a shame to base any proof on them. We
can prove anything by using them and they have caused so much
misery and created so many paradoxes. \ldots. Finally my eyes were
suddenly opened since, with the exception of the simplest cases,
for instance the geometric series, we hardly find, in mathematics,
any infinite series whose sum may be determined  in a rigorous
fashion, which means the most essential part  of mathematics has
no foundation. For the most part, it is true that the results are
correct, which is very strange. I am working to find out  why, a
very interesting problem.")

\bigskip

The author was struck by this citation when she first read it 26
years ago. At that time she knew nothing about divergent series. But this citation has followed her for her whole career.
In her point of view, this citation contains the past, present and
future of divergent series in mathematics.

\bigskip

\noindent {\bf The past:} As remarked by Abel, divergent series
occur very often in many natural problems of mathematics and
physics. Their use has permitted to do successfully a lot of
numerical calculations. One example of this is the computation by
Laplace of the secular perturbation of the orbit of the Earth
around the Sun due to the attraction of Jupiter. The calculations
of Laplace are verified experimentally, although the series he
used were divergent.

\bigskip

\noindent{\bf The present:} In the 20-th century divergent series
have occupied a marginal place. However it is during the same
period that we have learnt to justify rigorously their use,
answering at least partially Abel's question. In the context of
differential equations, the Borel summation, generalized by
\'Ecalle and others to multi-summability, gives good results (see
for instance \cite{Ba}, \cite{B}, \cite{E}).

\bigskip

\noindent{\bf The future:} Why do so many important problems of
mathematics lead to divergent series (see for instance \cite{I})?
What is the meaning of a series being divergent? Finding an answer to this question is a fascinating research field.

\bigskip

We will illustrate all this on the example of the Euler
differential  equation:
\begin{equation}
x^2y'+y=x. \label{Euler}\end{equation}

Since this is a short paper, the list of references is by no means
exhaustive.

\section{The past} The Euler
differential  equation \eqref{Euler} has the formal solution
\begin{equation} \hat{f}(x)=\sum_{n\geq0} (-1)^nn!
x^{n+1},\label{series}
\end{equation}
which is divergent for all nonzero values of $x$.

On the other hand, it is a linear differential equation whose
solution can be found by quadrature:
\begin{equation}
f(x)= e^{\frac1{x}}\ds\int_0^x\ds\frac{
e^{-\frac1{z}}}{z}dz.\label{solution}\end{equation} The integral
is convergent for $x>0$ and hence yields a solution of
\eqref{Euler}. We can rewrite this solution as
\begin{equation}
f(x)= \ds\int_0^x\ds\frac{
e^{\frac1{x}-\frac1{z}}}{z}dz,\label{solution2}\end{equation} in
which we make the change of coordinate $\frac{\zeta}{x}=\frac1{z}
-\frac1{x}$. This yields
\begin{equation}
f(x)= \int_0^{+\infty} \frac{
e^{-\frac{\zeta}{x}}}{1+\zeta}d\zeta.\label{solution3}\end{equation}
The integral is convergent for $x\geq0$. Hence, the solution $f(x)$ is
well defined for $x>0$ and, moreover, $\lim_{x\to 0^+}f(x)=0$.

\bigskip

What is now the link between the divergent power series
$\hat{f}(x)$ and the function $f(x)$? We will show that the
difference between $f(x)$ and a partial sum
\begin{equation} f_k(x)=\sum_{n=0}^{k-1}(-1)^n n!x^{n+1}\end{equation} of the series is smaller than the
first neglected term (this part has been inspired by \cite{R2}).
This is exactly as in the Leibniz criterion for alternating
series.

\bigskip

\begin{proposition}
 For any $x\geq0$
\begin{equation} |f(x)-f_k(x)|\leq k!
x^{k+1}.\label{approximation} \end{equation}
\end{proposition}
\noindent {\scshape Proof.} We use the following formula which is
easily checked
\begin{equation}
\frac1{1+\zeta} = \sum_{n=0}^{k-1} (-1)^n \zeta^n+
(-1)^{k}\frac{\zeta^k}{1+\zeta}.\label{geometric}\end{equation}
Then
\begin{equation}
\begin{array}{lll}
f(x)&=& \ds\int_0^{+\infty} e^{-\frac{\zeta}{x}}
\left(\sum_{n=0}^{k-1} (-1)^n\zeta^n +(-1)^k
\frac{\zeta^{k}}{1+\zeta}\right)d\zeta\\ &=&
\ds\sum_{n=0}^{k-1}\int_0^{+\infty} (-1)^n \zeta^n
e^{-\frac{\zeta}{x}}d\zeta + \int_0^{+\infty}
(-1)^k\frac{\zeta^{k}e^{-\frac{\zeta}{x}}}{1+\zeta}
d\zeta.\end{array}\end{equation}

Using that \begin{equation} \Gamma(n+1)= n! =
\int_0^\infty z^ne^{-z}dz, \label{gamma}\end{equation} which
implies \begin{equation}\int_0^\infty
\zeta^ne^{-\frac{\zeta}{x}}d\zeta= n!x^{n+1},
\label{gamma1}\end{equation}
 this yields
\begin{equation}
\begin{array}{lll}
\ds f(x) &=&\ds\sum_{n=0}^{k-1} (-1)^n n!x^{n+1}
+\ds\int_0^{+\infty}
(-1)^k\ds\frac{\zeta^ke^{-\frac{\zeta}{x}}}{1+\zeta} d\zeta\\
\ds&=& f_k(x) + R_k(x),\end{array}\end{equation} where $R_k(x)$,
the remainder, is the difference between $f(x)$ and the partial
sum $f_k(x)$ of the power series. The result follows since
\begin{equation}
\ds |R_k(x)|= \ds\int_0^\infty
\ds\frac{\zeta^ke^{-\frac{\zeta}{x}}}{1+\zeta} d\zeta\leq\ds\int_0^\infty \zeta^ke^{-\frac{\zeta}{x}}d\zeta= \ds
k!x^{k+1}.\label{remainder}\end{equation} \hfill $\Box$

\bigskip

The argument given here, which justifies rigorously the use of the
divergent series in numerical calculations, is very similar to the
argument made by Cauchy for the use of the Stirling series. In
particular, Cauchy used the formula \eqref{geometric}.

\bigskip

\noindent{\bf If we use the power series to approximate the
function $f(x)$, how good is the approximation?}

We encounter here the main difference between convergent and
divergent series. With convergent series, the more terms we take
in the partial sum, the better the approximation. This is no more the
case with divergent series. Indeed, if we take $x$ fixed in
\eqref{series} the general term tends to $\infty$. Hence, we are
better to take the partial sum for which the first neglected term
is minimum. In our example, $k!x^{k+1}$ is minimum when $x\sim \frac1{k}$, in which
case $\frac1{x}\sim k$. We use Stirling formula to approximate
$k!$ for $k$ large: \begin{equation} k!=k(k-1)!=k \Gamma(k)\sim
\sqrt{2\pi } k^{k+\frac12}e^{-k}.\end{equation} This gives us
\begin{equation}
|R_k(x)|\leq \sqrt{2\pi }\frac{e^{-k}}{\sqrt{k}}\sim \sqrt{2\pi x}
e^{-\frac1{x}},\end{equation} which is exponentially small for $x$
small. Not only have we bounded the error made when approximating
the function by the partial sum of the power series, but this
error is exponentially small for small $x$, i.e. very satisfactory
from the numerical point of view.

The phenomenon we have described here is not isolated and explains
the successes encountered when using divergent series in numerical
approximations.

\section{The present}

Looking at what we have done with the Euler equation, someone can
have the impression that we have cheated. Indeed, we have chosen a
linear differential equation, thus allowing us to construct by
quadrature a function which is a solution of the differential
equation. But what about the general case?

In general, once we have a formal solution by means of a power
series, we use a \emph{theory of resummation} for finding a function which is
a solution. In this paper, we will focus on the \emph{Borel method of
resummation}, also called 1-summability. We start with the
properties that an adequate theory of summability must satisfy:

\bigskip

\noindent{\bf Properties of a good method of resummation (see for
instance \cite{B}):} we consider a series $\sum_{n=0}^\infty a_n$,
to which we want to associate a number $S$, called its sum:

\begin{description}
\item{(1)} If $\sum_{n=0}^\infty a_n$ is convergent, then $S$ should be
the usual sum of the series.
\item{(2)} If $\sum_{n=0}^\infty a_n$ and $\sum_{n=0}^\infty b_n$ are
summable with respective sums $S$ and $S'$ then $\sum_{n=0}^\infty
(a_n+Cb_n) $ is summable with sum $S+CS'$.
\item{(3)} A series $\sum_{n=0}^\infty a_n$ is \emph{absolutely summable} if the series $\sum_{n=0}^\infty |a_n|$ is summable. 
\item{(4)} If $\sum_{n=0}^\infty a_n$ and $\sum_{n=0}^\infty b_n$ are
absolutely summable with respective sums $S$ and $S'$, then the product of the series, $\sum_{n=0}^\infty c_n$, where $c_n=\sum_{i=0}^n a_ib_{n-i}$, is absolutely summable with sum $SS'$.
\item{(5)} If $\sum_{n=0}^\infty a_n$ is summable with sum $S$,
then $\sum_{n=1}^\infty a_n$ is summable with sum $S-a_0$.
\item{(6)} (In the context of differential equations) If $\sum_{n=0}^\infty a_nx^n$ is summable with sum
$f(x)$, then $\sum_{n=1}^\infty n a_nx^{n-1}$ is summable with sum
$f'(x)$.
\end{description}

\bigskip

\noindent{\bf The Borel method of resummation for a divergent power series:} we
present it in a way which proves at the same time the property
(1): the idea is to take a convergent power series $\sum_{n=0}^\infty
a_nx^{n+1}$ and to write its sum $S(x)= \lim_{k\to \infty} \sum_{n=0}^k a_nx^{n+1}$
in a different way. For this purpose, we use \eqref{gamma} and we
write
\begin{equation}
\begin{array}{lll}
S(x)&=&\ds \sum_{n=0}^\infty a_nx^{n+1}  \\ &=&\ds\sum_{n=0}^\infty
\ds\frac{a_n x^{n+1}}{n!} \;n!\\ &=& \ds\sum_{n=0}^\infty\ds\frac{a_nx^{n+1}}{n!}
\ds\int_0^\infty z^n e^{-z} dz\\ &=& \ds\int_0^\infty
\left(\ds\sum_{n=0}^\infty\ds \frac{a_n (xz)^n}{n!} \right) e^{-z }
xdz\\
&=& \ds\int_0^\infty
\left(\ds\sum_{n=0}^\infty\ds \frac{a_n \zeta^n}{n!} \right) e^{-\frac{\zeta}{x} }
d\zeta.\end{array}
\end{equation}
Changing the order of the summation and the integral in the last line requires a proof, which we leave as an exercise. 

This suggests the following definition, which is well adapted to further analytic extension.

\begin{definition}\begin{enumerate}
\item A power series $\sum_{n=0}^\infty a_nx^{n+1}$ is 1-summable
(Borel-summable) in the direction $d$, where $d$ is a half-line
from the origin in the complex plane if 
\begin{itemize}
\item the series
$\sum_{n=0}^\infty a_n\frac{\zeta^n}{n!}$ has a nonzero radius of
convergence, and the sum of the series is an analytic function $g(\zeta)$ on the disk of convergence,
\item the function $g(\zeta)$ can be extended along the half line $d$,
\item and the integral
\begin{equation}
\int_dg(\zeta)
e^{-\frac{\zeta}{x}}d\zeta\end{equation} is convergent in some domain, with value
$S(x)$. \end{itemize} We call $S(x)$ the sum of the power series.\item A power series
$\sum_{n=0}^\infty a_nx^{n+1}$ is 1-summable if it is 1-summable
in all directions $d$, except a finite number of exceptional
directions.\end{enumerate}\end{definition}

 \noindent{\bf Example:}
In the case of the solution to Euler's differential equation we
have $a_n= (-1)^n n! $. Then
$$\sum_{n=0}^\infty \frac{a_n}{n!} \zeta^n= \frac1{1+\zeta} = g(\zeta).$$
Hence,
\begin{equation}S(x)=
\ds\int_0^\infty g(\zeta) e^{-\frac{\zeta}{z} } dz=
\ds\int_0^\infty\ds\frac{e^{-\frac{\zeta}{x}}}{1+\zeta}d\zeta,\label{solution4}\end{equation}
which is precisely \eqref{solution3}.
The solution of Euler's differential
equation is 1-summable in all directions except the direction
$\mathbb R^-$. The problem with the direction $\mathbb R^-$ comes
from the singularity at $\zeta=-1$  in \eqref{solution3} or
\eqref{solution4}.

\bigskip

\begin{figure}\begin{center}\includegraphics[width=6.5cm]{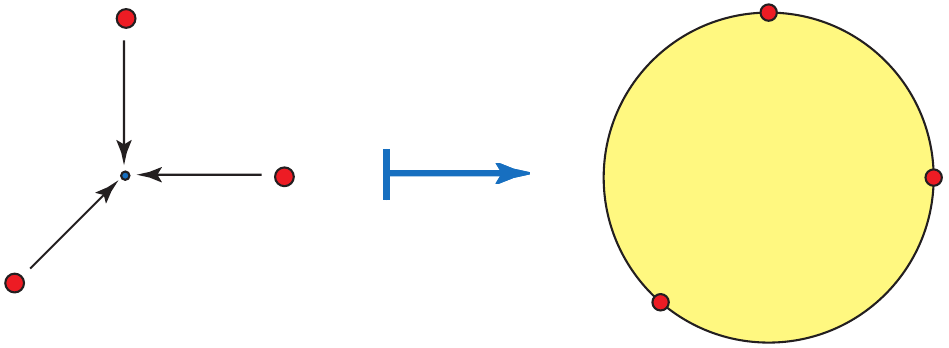}\caption{The Borel transforms sends the hidden singularities to a finite distance.}\label{fig:borel}\end{center}\end{figure}

\begin{figure}\begin{center}\subfigure[Convergent case]{\includegraphics[width=7cm]{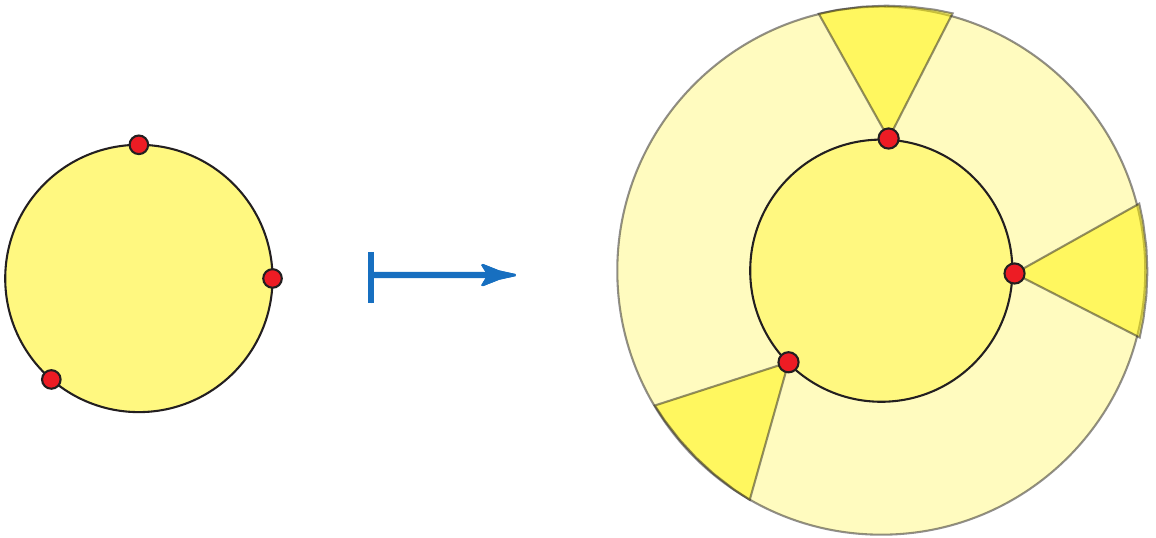}}\qquad\quad \subfigure[Divergent case]{\includegraphics[width=7cm]{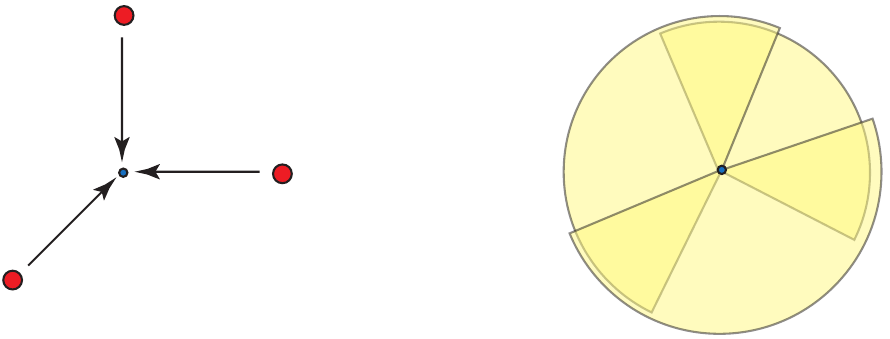}}\caption{We can extend the function around the singularities. }\label{extension}\end{center}\end{figure}

\noindent{\bf A very powerful insight coming from Borel.}  It is well known that if a power series $\sum_{n=0}^\infty
a_nx^{n+1}$ has a radius of convergence equal to $r$ and its sum
is a function $f(x)$ for $|x|<r$,  then the function $f(x)$ has a least one
singularity on the circle $|x|=r$. The idea of Borel is that a
divergent series is a series with radius of convergence $r=0$.
Hence, we have at least one singularity hidden in some direction:
for Euler's differential equation, this is the direction $\mathbb R^-$. The
operation \begin{equation} \hat{f}(x)= \sum_{n=0}^\infty
a_nx^{n+1} \qquad \mapsto\qquad \mathcal{B}(\hat{f})=
\sum_{n=0}^\infty
\frac{a_n\zeta^n}{n!},\label{Borel_trans}\end{equation} sends the
singularities at a finite distance, as if we had blown the disk of convergence to bring it from radius $0$ to a finite radius $r$ (see Figure~\ref{fig:borel}). $\mathcal{B}(\hat{f})$ is called
the {\it Borel transform} of $\hat{f}$. For a convergent power series on a disk, we can extend the function outside the disk of convergence, as long as we avoid the singularities. In particular, we can turn around the singularities (see Figure~\ref{extension}(a)). Generically, the extension is ramified at the singularities. If we have only a finite number of singularities and we now shrink the radius of the disk to $0$, then we get a function defined on a finite number of sectors covering a pointed neighborhood of the origin (see Figure~\ref{extension}(b)), which is what is provided by the Borel summability.

\bigskip

A theory of resummation is useful if there are theorems associated
to it. For instance, for the Borel summation, let us cite this
theorem of Borel, which has been followed by many others theorems of the same type.

\begin{theorem} \cite{B} We consider an algebraic differential
equation \begin{equation}F(x,y,y',\dots
y^{(m)})=0,\label{algebraic}\end{equation} where $F$ is a
multivariate polynomial. If $\hat{f}(x)=\sum_{n=0}^\infty a_nx^n$
is a formal solution of \eqref{algebraic} and $\hat{f}$ is
absolutely Borel-summable with Borel sum $f(x)$, then $f(x)$ is
solution of the differential equation \eqref{algebraic} and has
the asymptotic expansion $\hat{f}(x)=\sum_{n=0}^\infty
a_nx^n$.\end{theorem}

\bigskip

\begin{remark} \begin{description}
\item{(1)} The sums of a 1-summable series in the different
directions $d$ give functions which are analytic extensions one of
the other as long as we move the line $d$ continuously through
directions in which the series is 1-summable. This yields a
function defined on a sector with vertex at the origin.
 More details can be found for instance in
\cite{Ba}, \cite{R1} and \cite{R2}.
\item{(2)} The Borel sum of a divergent power series can never be uniform in
a punctured neighborhood of the origin. It is necessarily
ramified. This is what is known in the literature as the {\it
Stokes phenomenon}.
\end{description}
\end{remark}

The 1-summability and its extensions have been extensively studied
during the 20-th century.  Extensions of the notion of
1-summability have been obtained by allowing ramifications
$x=\xi^m$. Then 1-summability in $x$ corresponds to
$m$-summability in $\xi$. The notion of multi-summability has also
been introduced: a series $\hat{f}$ is multi-summable if it is a
finite sum $\hat{f}=\hat{f}_1+\dots +\hat{f}_n$, each $\hat{f}_i$
being $m_i$-summable. Explicit criteria allow deciding a priori
that some systems of differential equations have multi-summable
solutions in the neighborhood of certain singular points (see for
instance \cite{Ba}).

\section{The future}\label{sec:future}
Let us now look at a generalized Euler differential equation
\begin{equation} x^2y'+y=g(x).\label{Euler_gen}
\end{equation} For almost all functions $g(x)$ analytic in the
neighborhood of the origin and such that $g(0)=0$  the formal
solution of \eqref{Euler_gen} vanishing at the origin is given by
a divergent series.   Only in very special cases the equation
\eqref{Euler_gen} has an analytic solution at the origin. For instance, $f(x)=x$ is the analytic solution of
\begin{equation} x^2y'+y=x+x^2.\label{Euler2}\end{equation}

\bigskip

What is the difference between the equations \eqref{Euler} and
\eqref{Euler2}? To understand we apply successively the two
steps:\begin{itemize} \item Complexification: we allow $x\in
\mathbb C$;
\item Unfolding: $x=0$ is a double singular point of each equation. Hence, we introduce a parameter so as to split the double singular point into two simple singular points.
\end{itemize}

\bigskip

\begin{figure}[h]
\begin{center}
\subfigure[ Equation \eqref{Euler}]
    {\includegraphics[angle=0,width=6cm]{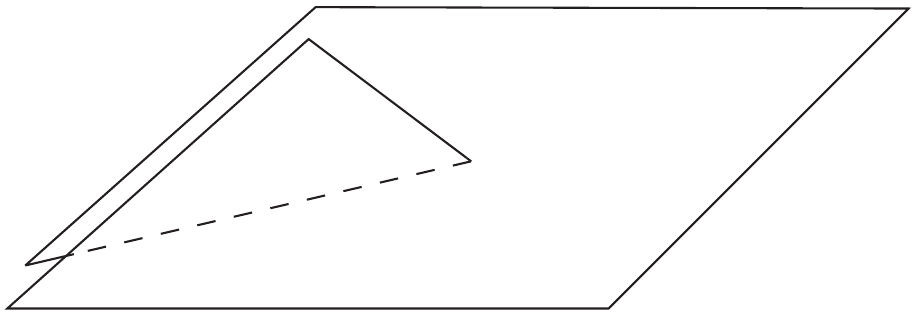}}
     \qquad\qquad
\subfigure[ Equation \eqref{Euler2}]
    {\includegraphics[angle=0,width=6cm]{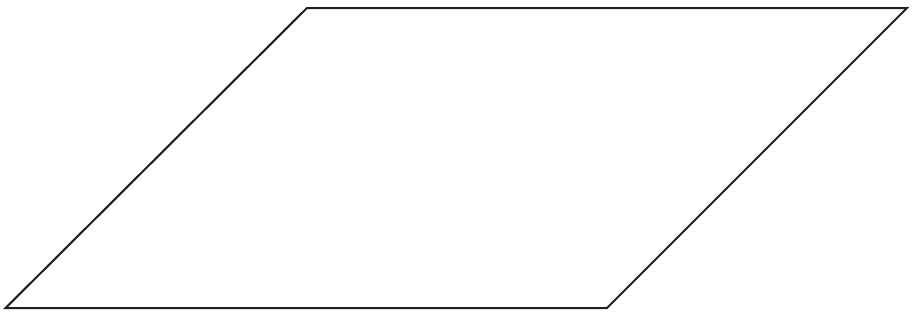}}
\caption{The domain of $f(x)$} \label{fig1}
\end{center}
\end{figure}

\noindent{\bf Complexification:} let us consider a function $f(x)$
which is the Borel sum of a solution of \eqref{Euler_gen}, and its analytic extension when we turn around the origin.
The function $f(x)$ is (Figure~\ref{fig1}):\begin{itemize}\item
uniform for \eqref{Euler2};
\item ramified for \eqref{Euler}, and generically ramified for a
solution of \eqref{Euler_gen}. The two branches differ by a
multiple of $e^{\frac1{x}}$ (which is a solution of
the homogeneous equation).
\end{itemize}
\bigskip

\noindent{\bf How to understand that generically we should have
ramification?}
 To answer, we unfold and embed \eqref{Euler_gen} into
 \begin{equation}
 (x^2-\eps)y'+y=g_\eps(x).\label{Euler_unfold}\end{equation}
  We limit here our discussion to
 $\eps>0$, although all complex values of $\eps$ are of interest.
The solutions of the homogeneous equation associated to
\eqref{Euler_unfold} are given by
$C\left(\frac{x+\sqrt{\eps}}{x-\sqrt{\eps}}\right)^{\frac1{2\sqrt{\eps}}}$.
We localize at the singular points using the changes of coordinate $X_\pm= x\mp \sqrt{\eps}$, thus transforming the equation to the form
$$X_\pm y'+ \frac{y}{X_\pm \pm 2\sqrt{\eps}}= k_\eps^\pm(X_\pm).$$

\medskip A particular solution $y=h_1(x)$ near $x_1=\sqrt{\eps}$ is analytic. Indeed, taking $X_+=x-\sqrt{\eps}$, it is given by 
 $$h_1(x)= \left(\ds\frac{X_++2\sqrt{\eps}}{X_+}\right)^{\frac1{2\sqrt{\eps}}}\int_0^{X_+}k_\eps^+(z) z^{\frac1{2\sqrt{\eps}}-1}(z+2\sqrt{\eps})^{-\frac1{2\sqrt{\eps}}}dz.$$
 The general solution is given by 
 $$y(x) = h_1(x) + C_1 \left(\frac{x+\sqrt{\eps}}{x-\sqrt{\eps}}\right)^{\frac1{2\sqrt{\eps}}}.$$
 Hence, the solution $h_1(x)$ (corresponding to $C_1=0$) is the unique solution which is analytic and bounded at $x=\sqrt{\eps}$.

 \medskip Similarly, if we let $X_-=x+\sqrt{\eps}$,  a particular solution $y=h_2(x)$  near $x_2=-\sqrt{\eps}$ is given by 
 $$h_2(x)= \left(\ds\frac{X_-}{X_--2\sqrt{\eps}}\right)^{\frac1{2\sqrt{\eps}}}\int_0^{X_-}k_\eps^-(z) z^{-\frac1{2\sqrt{\eps}}-1}(z-2\sqrt{\eps})^{\frac1{2\sqrt{\eps}}}dz.$$
 The general solution is also given by 
 $$y(x) = h_2(x) + C_2 \left(\frac{x+\sqrt{\eps}}{x-\sqrt{\eps}}\right)^{\frac1{2\sqrt{\eps}}}.$$

We now have two cases:
\begin{description}\item{(1)} If $\frac{1}{2\sqrt{\eps}}\notin
\mathbb N$, then $h_2$ is analytic near $x=-\sqrt{\eps}$. Since $
C_2\left(\frac{x+\sqrt{\eps}}{
 x-\sqrt{\eps}}\right)^{\frac1{2\sqrt{\eps}}}$ is ramified at $x=-\sqrt{\eps}$ for
$C_2\neq0$, all solutions but $h_2$ are ramified. It is of
course exceptional that the extension of $h_1$ at
$-\sqrt{\eps}$ be exactly the solution $h_2$ and, generically,
we should expect that the analytic extension of $h_1(x)$ is $h_2(x)+
C_2\left(\ds\frac{x+\sqrt{\eps}}{
 x-\sqrt{\eps}}\right)^{\frac1{2\sqrt{\eps}}}$ with $C_2\neq0$. Hence the extension of $h_1$ should be ramified at
$x=-\sqrt{\eps}$. If this ramification holds till the limit
$\eps=0$ we get Figure~1(a).
\item{(2)} Let us now consider the case $\frac{1}{2\sqrt{\eps}}=n\in
\mathbb N$. We must again divide the discussion in two cases:
\begin{description}\item{(i)} In the generic case, $h_2$ is ramified: it
contains one term of the form
$(x+\sqrt{\eps})^n\ln(x+\sqrt{\eps})$.
\item{(ii)} In the exceptional case, $h_2$ is analytic, and so are all solutions through
$x=-\sqrt{\eps}$. Hence, it is impossible for $h_1$ to be
ramified at $x=-\sqrt{\eps}$. This case is excluded in the
unfolding of an equation \eqref{Euler_gen} whose solution is
ramified.
\end{description}
\end{description}

 Let us now summarize the very interesting phenomenon we have discovered: if the
formal solution of \eqref{Euler_gen} is divergent, then in the
unfolding:
\begin{itemize}
\item Necessarily $h_1(x)$ is ramified at $-\sqrt{\eps}$: the
divergence means a form of \emph{incompatibility} between the special local
solutions at the two singular points $\pm\sqrt{\eps}$, which remains
until the limit at $\eps=0$. \item {\bf Parametric resurgence phenomenon:} for
sequences of values of the parameter $\eps$ converging to $\eps=0$
(here $\frac1{2\sqrt{\eps}}\in \mathbb N$), the pathology of the
system is located exactly at one of the singular points. Indeed,
here, the only way for the system to have an incompatibility is
that one of
 the singular points be pathologic.
\end{itemize}
\bigskip

We have understood why divergence is the rule, and convergence the
exception.  

\bigskip

\begin{figure}[ht!]\begin{center} \includegraphics[width=5cm]{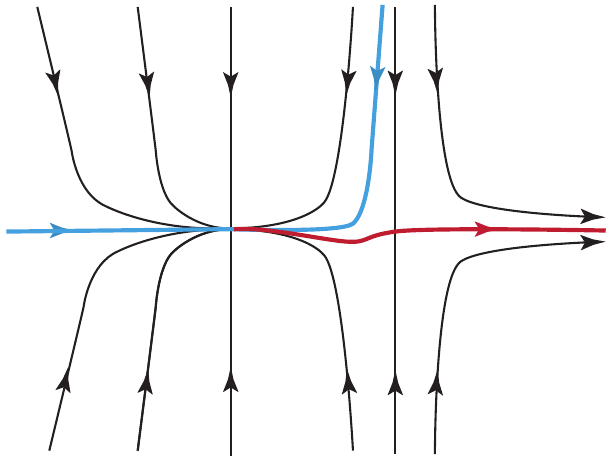}
\caption{The analytic invariant manifolds of the saddle and the node in \eqref{saddle-node_unfold} do not match.} \label{nomatch} \end{center} \end{figure}
\begin{figure}
\begin{center}
\subfigure[$\eps=0$]
{\includegraphics[height=4cm]{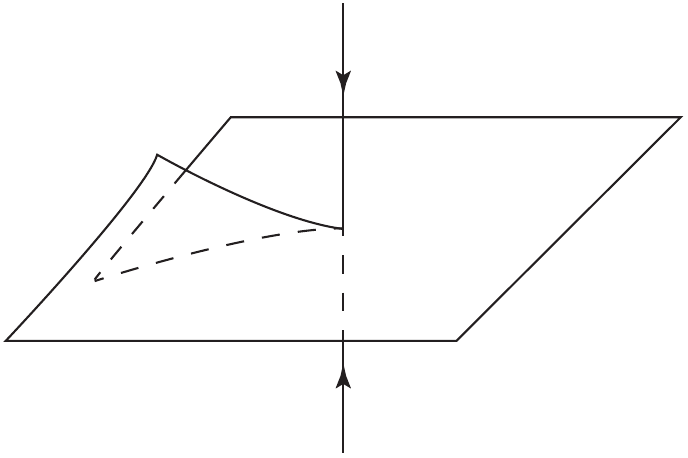}}\qquad\quad
\subfigure[$\eps>0$]{\includegraphics[height=4cm]{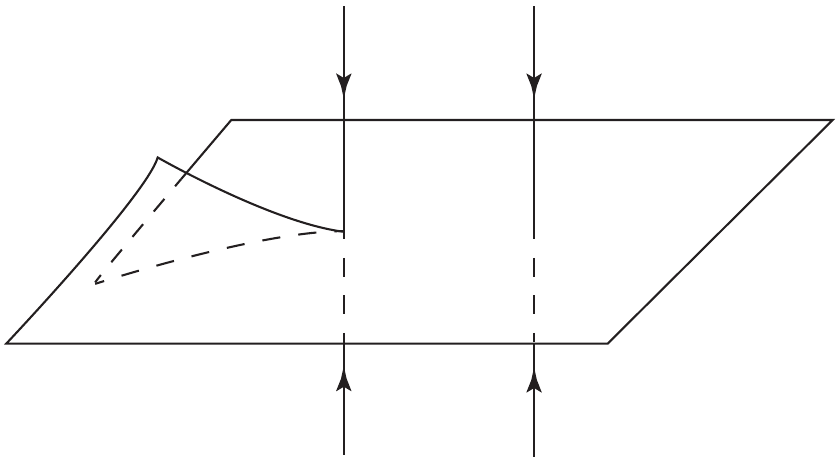}}\caption{A non analytic center manifold of a saddle-node \eqref{saddle-node_unfold} and its unfolding: the $x$-coordinate is drawn in $\C$.}\label{center_manifold}\end{center}\end{figure}

\noindent{\bf Application.} The generalized Euler differential equation \eqref{Euler_gen} is the equation for the center manifold of the saddle-node in the two-dimensional analytic  vector field
\begin{align}\begin{split} 
\dot x&=x^2,\\
\dot y &= - y +g(x),\end{split}\end{align} 
with $g(0)=0$. More generally, if we have a saddle-node in a (real) two-dimensional analytic  vector field
\begin{align}\begin{split} 
\dot x&=x^2,\\
\dot y &= - y +O(x^2) + O(y^2).\end{split}\end{align} 
it is known that, generically, the center-manifold is only $\C^\infty$. Now, we understand why. 
Indeed, the unfolded system \begin{align}\begin{split} 
\dot x&=x^2-\eps,\\
\dot y &= - y +O(x^2-\eps) + O(y^2),\end{split}\label{saddle-node_unfold}\end{align} 
has two singular points at $(\pm\sqrt{\eps},0)$ with ratio of eigenvalues $\lambda_\pm(\eps)=\pm 2\sqrt{\eps}\left(1+ O(\sqrt{\eps})\right)$.  When $\eps>0$, the point $(\sqrt{\eps},0)$ is a saddle, and the point $(-\sqrt{\eps},0)$ is a node. The saddle has an analytic unstable manifold. For $\eps\notin \frac1{\N}$, then the node is linearizable. In the linearizing coordinates $(X,Y)$, the trajectories have the form $Y=CX^{\frac1{2\sqrt{\eps}}}$ and $X=0$.  Hence, the node has exactly two invariant analytic manifolds $X=0$ and $Y=0$. There is no reason why $Y=0$ would match with the unstable manifold of the saddle (see Figures~\ref{nomatch} and \ref{center_manifold}). If it does not match, then, when taking complex coordinates, the unstable manifold of the saddle is ramified at the node. If the mismatch is up to the limit, then we expect the center manifold to be not analytic at the limit, and only $1$-summable.  

When $\eps\in \frac1{\N}$, then we have the parametric resurgence phenomenon. If the saddle-node at the limit has no analytic center manifold, then the node is resonant with normal form
\begin{align}\begin{split} 
\dot X&=-\frac1n X,\\
\dot Y &= - Y + a_nX^n,\end{split}\end{align} 
with $a_n\neq0$. 

\section{Conclusion} The phenomena described in Section~\ref{sec:future} are much more general than the context
of \eqref{Euler_gen} discussed here, and are explored by the author and her collaborators in different contexts: parabolic points of 1-dimensional diffeomorphisms, saddle-nodes and resonant saddles of two-dimensional vector fields, Hopf bifurcation in two-dimensional vector fields, non resonant irregular singular points.  
A common thread is that $1$-summability occurs when two equilibrium points of a dynamical system coallesce in a double point, or more generally two special \lq\lq  invariant objects\rq\rq\ (for instance a singular point and a limit cycle) coallesce. 
There are some rigid dynamics attached to each simple equilibrium (object), and these rigid dynamics do not match well until the limit when the two equilibria (or objects) coallesce. 
More generally, $k$-summability occurs when $k+1$ equilibrium points (or objects) merge together in a generic way. Multi-summability can occur when the merging is non generic. 
\bigskip

Divergent series also occur in the  phenomena involving \lq\lq
small denominators\rq\rq. A source of divergence in this case is the accumulation of particular solutions. For instance, in the case of a fixed point of a
germ of analytic diffeomorphism $f:(\mathbb C,0)\rightarrow\mathbb
(\mathbb C,0)$, $f(z) =e^{2\pi i \alpha}z + o(z)$, with $\alpha\in
\mathbb R\setminus\mathbb Q$, the divergence of the linearizing
series could come from the accumulation of periodic points around the fixed point at the origin. The study of this phenomenon is part of the work for which Jean-Christophe Yoccoz received
the Fields medal in 1994.
 
\bigskip
Divergent series occur generically in many situations within
differential equations and dynamical systems. The divergence of the series
carries a lot of geometric information on the solutions. For instance, if the formal power series solution of
the Euler equation \eqref{Euler} had been convergent, its sum
could not have been ramified. 
\bigskip

\noindent{\bf The future:} \begin{itemize} \item A large program
of research will consist in learning  to \lq\lq read\rq\rq\ all the rich
behaviour of functions defined by divergent series.
\item Shouldn't divergent series occupy a more important place in
mathematics?
\end{itemize}

\end{document}